\documentclass[12pt]{article}
\usepackage{amsmath}
\usepackage{amssymb}
\newtheorem{theorem}{Theorem}[section]
\newtheorem{definition}[theorem]{Definition}

\newtheorem{lemma}[theorem]{Lemma}
\numberwithin{equation}{section}

\newcommand{\gd}{\dot{\frak g}}
\newcommand{\C}{\mathbb C}
\newcommand{\Z}{\mathbb Z}

\newcommand{\dg}{\dot{\frak g}}
\newcommand{\K}{{\cal K}}
\newcommand{\R}{{\cal R}}
\newcommand{\dd}{\partial}
\newcommand{\D}{{\cal D}}
\newcommand{\m}{{\bf m}}
\newcommand{\rb}{{\bf r}}

\newcommand{\g}{\frak g}
\newcommand{\gl}{\widehat{gl}_N}
\newcommand{\slN}{\widehat{sl}_N}
\newcommand{\h}{{Hei}}

\newcommand{\V}{{\cal V}}

\newcommand{\ma}{{\m \cdot {\bf a}}}
\newcommand{\ra}{{\rb \cdot {\bf a}}}
\newcommand{\mra}{{(\m+\rb) \cdot {\bf a}}}
\newcommand{\bs}{b_s}

\newcommand{\bq}{b_q}

\newcommand{\wdg}{\widehat \dg}

\newcommand{\dV}{\dot V}
\newcommand{\Lor}{Lat}
\newcommand{\N}{\mathbb N}

\begin{document}

\title{Vertex operator algebras and the representation theory of toroidal algebras}

\author{Stephen Berman\thanks{Research supported by the Natural Sciences and
Engineering Research Council of Canada}\\
Department of Mathematics and Statistics\\University of Saskatchewan, 
Saskatoon, Canada\and Yuly
Billig${}^*$\\
School of Mathematics and Statistics\\ Carleton University, Ottawa, Canada \and Jacek
Szmigielski${}^*$\\ Department of Mathematics and 
Statistics \\University of
Saskatchewan\\ Saskatoon, Canada}

 \maketitle

{\abstract  An explicit vertex operator algebra construction is given of a class of irreducible 
modules for toroidal Lie algebras. }

\medskip

\medskip

\section{Introduction} This paper is about representations of toroidal Lie algebras and vertex operator algebras  which are naturally associated to them. In fact, we find the picture is much like the case of the affine Kac-Moody Lie algebras, where, for quite some time now, 
VOAs have played an important role in the study of their representations. Indeed, this link with VOAs and affine Kac-Moody Lie algebras has been one of the driving forces for both the representation theory of these algebras as well as for  the development of the  theory of VOAs. Now toroidal Lie algebras are multi-variable generalizations of the affine algebras and so it is extremely natural to consider VOAs associated to their representation theory. In fact, roughly speaking, one obtains a toroidal Lie algebra associated to any finite dimensional simple Lie algebra over $\C$, call it $\gd$, as follows. One tensors $\gd$ with the Laurent polynomials in several, say $N+1$ where $N\geq 1,$ variables to get a ``multi-loop'' algebra and then takes the universal central extension of this algebra. Next, just as in the affine case, one adds some derivations. The result is a toroidal algebra. Here however, the reader should understand that already this process is quite a bit more complicated than the affine case due to the fact the universal central extension will always be infinite dimensional as $N+1 \geq 2$ and hence the derivation algebra we need to add will also be infinite dimensional. Letting $\R$ denote this ring of Laurent polynomials one knows,\cite{Kas}, that the space of 1-forms modulo exact forms, we will denote this space by $\K$, needs to be added to $\gd \otimes \R$ to get the universal central extension and that the algebra $\D^*$ which we need to add is a subalgebra of $\text{Der}(\R).$ Thus, our toroidal Lie algebra looks like
\begin{equation}
\gd \otimes \R \oplus \K \oplus \D^*.
\end{equation}

In the first section of this paper we will define all of this rigorously. Here, we want to mention that, in fact, there is a possible cocycle 
\begin{equation}
\tau:\D^* \times \D^* \rightarrow \K,
\end{equation}
which enters into the definition of the toroidal algebra. Thus, we denote the resulting toroidal algebra by
\begin{equation}
\g_\tau =\gd \otimes \R \oplus \K \oplus \D^*.
\end{equation}
In our set-up these possible cocycles come from a $2$-dimensional space and hence parametrize the class of toroidal algebras
 which we will work with here. One of these cocycle was first discovered in \cite{RM} where the authors gave a representation, via vertex operators, for one of the algebras, $\g_\tau$. This cocycle was produced by the representation itself and the authors of this work did the appropriate computations necessary to see exactly what this cocycle was when a particular algebra of derivations, $\D^*$, was used. This has influenced the entire theory as right from its start  both vertex operators and cocycles were present.

In the work \cite{B} the author gave what could be called the principal realization by vertex operators for this same toroidal algebra with the cocycle discovered in \cite {RM}. This time the associated representation was irreducible and 
the space of the representation, the so-called Fock space, was just a ring of Laurent polynomials in $N$ variables tensored with a symmetric algebra in infinitely many variables. The work \cite{BB}, building on \cite{B} as well as on \cite{R,La}, produced infinite families of irreducible representations for all of the algebras $\g_\tau$. Moreover, these representations have finite dimensional weight spaces relative to a natural Cartan subalgebra of the toroidal algebra $\g$, and besides this, they are as faithful as is possible. What this means is the following. The toroidal algebra has an
 $N+1$ dimensional centre (most of the elements from $\K$ do not remain central when $\D^*$ is added) and so in any representation which is irreducible and has finite dimensional weight spaces one must have a kernel which is {\bf {at least}} $N$ dimensional by Schur's Lemma. The representations constructed in \cite{BB} have kernels which are exactly $N$ dimensional. Thus, they are as faithful as they possibly can be. We proceed to describe them in a little more detail here.

   Let $\V=\K \oplus \D^*$ and let $W$ be any finite dimensional $gl_N$-module. Next, let $T(W)=\R_N \otimes W$ where $\R_N$ is the ring of Laurent polynomials in $N$ variables. Then, as seen from \cite {R,La}, $T(W)$ becomes a module for $\text{Der}(\R_N)$. This so called tensor module corresponds geometrically to the action of vector fields on tensor fields by the Lie derivative. Now $\V$ has a natural $\Z$ gradation which corresponds to the first variable, $t_0$ in the natural $\Z^{N+1}$-gradation of $\V$ and hence this allows us to write
\begin{equation}
\V=\V_-\oplus \V_0 \oplus \V_+
\end{equation}
where $\V_0$ is all elements of degree $0$ and $\V_+$ (respectively $\V_-$) is all elements of positive (respectively negative) degree. Now it is not hard to see that $\V_0$ acts on $T(W)$ and this introduces  a scalar $c$, called the central charge, into this construction. Moreover, it is true that $T(W)$ is an irreducible $\V_0$ module provided only that $W$ is an irreducible $gl_N$ module and $c \neq 0$, see \cite{BB} Proposition 2.25. Just as in the case of affine Kac-Moody Lie algebras one forms the $\V$ module
\begin{equation}
\overline U(W,c)={\text{Ind}}_{\V_+\oplus\V_0}^{\V}T(W)
 \end{equation} 
(these modules were called $M(W,c)$ in \cite{BB}).

It is not hard to see that this $\V$ module has a unique submodule intersecting the top $T(W)$ trivially and if $W$ is chosen to be irreducible then the corresponding quotient module, denoted $L(W,c)$, is irreducible and has finite dimensional weight spaces (see \cite{BB},Theorem 4.9). The finite dimensionality is quite non-trivial and is seen because of the notion of polynomial multiplication which is introduced in \cite{BB}.

   In the final stage of this construction one needs to take the special derivation algebra 
\begin{equation}
\D^*=\{ \sum_{p=1}^N f_p d_p | f_p \in \R \},
\end{equation} 
where $d_p$ denotes the degree derivation of $\R$ corresponding to the $p$th variable $t_p$ but where any of the cocycles $\tau$ are allowed. Then we tensor one of the modules obtained above for $\V$ with a highest weight module for the non-twisted affine Lie algebra associated to $\gd$. Here both modules must have the same central charges. We then obtain modules for the toroidal algebra $\g_\tau$ which again, if our initial data satisfies that $W$ is irreducible,  $c \neq 0$ and the module for the affine algebra is irreducible, then yields an irreducible module for our toroidal algebra with finite dimensional weight spaces and which is as faithful as possible in the sense already described.

   In many ways the above picture resembles the representation theory of highest weight modules for the affine Kac-Moody Lie algebras and hence, it became clear that one should look for VOAs here as well. Looking at the representation given in \cite{B} one certainly expects an affine VOA, as well as a sub-VOA of a a hyperbolic lattice VOA to play a role here. There is also the possibility of duplicating the theory of affine VOAs as presented in \cite{FZ,Li} in this toroidal setting and in fact , in our first attempt we did exactly this. However, this approach only identifies some of the above modules as VOAs and so does not at all address the possibility of giving realizations for these modules using VOA theory. This is in fact what we accomplish in this work. We are able to give  VOAs $V$ which appear as a tensor product of three fairly well-known VOAs and then we show that this VOA is isomorphic to one of the above irreducible toroidal modules. Moreover, just as in the affine case, we are able to show that the toroidal modules introduced in \cite{BB}, will in fact be modules for this VOA. 
Since the characters of the three pieces in our realizations are already known,
it becomes clear that the characters of the modules from \cite{BB} are
just products of these known characters.

  We next describe the three types of VOAs which will make up $V.$ Here we only describe these VOAs roughly but later in the body of the paper we shall develop precise notation for all of them.
 Let $\widehat \gd$ be the untwisted affine algebra, with central element denoted $k_0$, built on $\gd$ and let $V_{\widehat \gd}$ be the Verma VOA associated to $\widehat \gd$ ($V_{\widehat \gd}$ is a generalized Verma module for $\gd$)
with central charge $c$. Here we assume $c\neq -h^\vee$ where $h^\vee$
is the dual Coxeter number of $\gd$. Then one knows that the 
quotient module $L_{\widehat \gd}$ is an irreducible highest weight module of $\widehat \gd$. Next, let $\Lor$ be the hyperbolic lattice with basis $\{a_p,b_p |1 \leq p \leq N \}$ where the $a_p$'s and $b_p$'s span isotropic subspaces and satisfy $(a_p,b_q)=\delta_{p,q}.$ We let $V_{\Lor}$ be the lattice VOA associated to this even lattice. Then we have 
\begin{equation}
V_{\Lor}=S_H \otimes \C[\Lor],
\end{equation}
where $S_H$ is a symmetric algebra and $\C[\Lor]$ is a twisted group algebra of the group algebra of $\Lor.$ We then let
\begin{equation}
V_{\Lor}^+=S_H \otimes \C[\Lor^+]
\end{equation}
be the sub-VOA associated to the isotropic subspace, $\Lor^+$, of $\Lor$ spanned by the elements $a_p,1 \leq p \leq N.$
Our third VOA is associated to the Lie algebras ${\gl}$ which we take to be the direct sum of $\slN$  with the central element $C_1$ and a Heisenberg Lie algebra, $\h$ 
 generated by the identity matrix and the central element $C_2$ so that
\begin{equation}
\h= I\otimes \C[t_0,t_0^{-1}] \oplus \C C_2.
\end{equation}
Then the VOA associated to this  algebra $\gl$ is taken to be the tensor product of the Verma VOA associated to the affine Kac-Moody Lie algebra ${\widehat 
{sl}_N}$ with the simple Heisenberg VOA associated to the Heisenberg algebra $\h$ above. We denote the first of these by $V_{{\widehat {sl}_N}}$ and the second one by $L_{\h}$. Thus, we are letting
\begin{equation}
V_{{\widehat gl_N}}=V_{{\widehat sl_N}} \otimes L_\h,
\end{equation}
be the VOA associated to the algebra ${\widehat {gl}_N}$. Then our VOA $V$ is given as
\begin{equation}
V= V_{\widehat \gd} \otimes V_{\Lor}^+ \otimes V_{{\widehat gl_N}}.
\end{equation}
Actually, because of the various choices for central charges, namely $c,c_1,c_2$ from above, we obtain a family of VOAs. Moreover, we will show that $c_1,c_2$ correspond to different choices of the cocycle $\tau$, while $c$ 
determines a cocycle associated with the action of the Virasoro algebra on $\V$.  All of this is explained later in the paper when the proper notation and results have been established.

   The proof of our result depends on carrying out certain computations in the VOA $V$. What we need to do is to define the proper vertex operators, $Y(v,z)$ for certain $v \in V$ which will be the images of the specified generating fields of the toroidal algebra, and then to show they satisfy the same commutation relations. This involves some fairly delicate computations, but these computations are straightforward enough for those well versed in VOAs. After this is done, the way is clear to show that modules for the VOA $V$ correspond to modules for the toroidal algebra $\g_\tau$ and to then explain how the fairly well-known irreducible modules for the VOAs involved provide a model for the modules from \cite{BB}. It is really only in this last part that a good knowledge of \cite{BB} is necessary. We have recalled here, in this Introduction, much from that work only to show what our motivation for the present work is.

   In the second section of this work we give the definition of the toroidal algebras as well as the cocycles $\tau$. We then present the generating fields which we use and display their commutation relations. One of these fields $K_0({\bf m},z)$ is special and it has conformal weight different from that of the other central fields.  We close this section with a simple Lemma which allows us to conclude certain commutation relations hold once we know that some others do. This lets us cut down on the amount of computation we must do. 

In the next section we begin by recalling the definition of a VOA. Here, we are dropping the requirement that the grading consists of {\bf finite dimensional} homogeneous spaces but we do require all other axioms. But note that when we identify the modules from \cite{BB} we will know that then the homogeneous pieces involved do have a $\Z^N$ gradation which has {\bf finite dimensional} homogeneous spaces. We then go on to discuss the particular types of VOAs we need. As above, these include affine VOAs, lattice VOAs, and Heisenberg VOAs. Most of this material is by now quite well known, so we are brief, and only develop the material to the point where we have enough notation and information to accomplish the computations which we do in the following section.

   The fourth section is the heart of the paper and here is where we present the details of our computation. We must check the various commutation relations case by case. We do this in such a way so as to give the reader a feeling for working with the various VOAs involved and to show how the various pieces from the tensor product of our three VOAs fit together to produce the toroidal algebra. Almost everything important that is taking place here is involved with {\bf negative} powers of $z$ so we often just work with the ``anti-holomorphic'' part of our series. We indicate this when we do this. After accomplishing this we go on, in the final section , to tie up our results to the modules from \cite{BB}.  We point out how irreducible modules for our VOA $V$ coincide with irreducible modules for $\g_\tau.$

   We close this introduction with several remarks. First, we want to mention that there is another very nice and quite different work by K. Saito and D. Yoshii on toroidal algebras and VOA's, \cite {SY}. They work with a lattice VOA attached to the full hyperbolic lattice and use the realization given in \cite{RM,MEY}. They then look at a certain Lie algebra attached to this VOA, namely the VOA factored by $L(-1)$ of it, and find the toroidal algebra is isomorphic to a subalgebra of this Lie algebra. They do not use any derivations of the toroidal algebra and so just work with the algebra $\gd \otimes \R \oplus \K$. One feature of their work is to give different descriptions of the same toroidal Lie algebra. 
The second remark we want to make here is that the part of our VOA $V$ associated to $\gl$ may seem a bit mysterious to some. It was in the insightful work of T. Larsson \cite {La,L} that this part seems to have been first understood. This insight of Larsson was built into \cite{BB} since, at the very beginning of the construction there, one begins with a finite dimensional
$gl_N$ module. This then naturally leads to the VOA $V_{{\widehat {gl}_N}}$. 

Regarding the references to a fast growing literature on VOAs we mention only 
those which we used directly in this work, to wit, the references \cite {FHL,FLM,K,Li,MN} were all quite helpful to us at various stages of this work and served as our basic sources for information on VOAs. 
In addition, the course given at the Fields Institute by Chongying Dong during the fall term 
of 2000 was extremely valuable to all three of us.  

   While working on this project the authors greatly benefited from many very helpful discussions with both Chongying Dong and Yi-Zhi Huang. This project certainly benefited from their help, encouragement, and insights, and the authors want to take this opportunity to offer thanks to them for this.

\section{Toroidal algebras}\label{Toroidalalg}

   For information and notation on toroidal algebras we will follow the work \cite {BB}. However we recall what we need here.

Let $N \geq 1$ be an integer and let ${\cal R}={\mathbb C}[t_0^{\pm 1},t_1^{\pm 1},...,t_N^{\pm 1}]$
  be the ring of Laurent polynomials in $N+1$ variables over the complex field and let $\gd$ be a finite dimensional Lie algebra over $\C.$ We let ${\cal K}= \Omega_{\cal R} ^1/d {\cal R}$ be the space of 1-forms modulo the exact forms. For notation we let $fdg$ denote the element of $\K$ corresponding to the pair of elements $f,g $ from $\R.$ We  let $k_i =t_i^{-1}d t_i$ .Then we have that  $\cal K$ is spanned by elements of the form ${t^{\bf m}k_i}$ where ${\bf m} \in \Z^{N+1}$,  $ 0 \leq i \leq N$.
 Here we are using for notation that  $t^{\bf m}=t_0^{m_0}t_1^{m_1} \dots t_N^{m_N}$ where
${\bf m} = (m_0, m_1, ...,m_N) \in \Z^{N+1}.$ The  relations staisfied by these elements are that
\begin{equation}\label {cent}
\sum_{p=0}^{ N} m_p t^{\bf m}k_p =0 , \quad {\bf m} \in \Z^{N+1}.
\end{equation}

 The space $\K$ is important for giving the universal central extension of the Lie algebra $\dg \otimes \R$, 
where the bracket in $\dg \otimes \R \oplus \K$ is given by
\begin{equation}\label {bracket}
[g_1 \otimes f_1(t), g_2\otimes f_2(t)] = [g_1, g_2] \otimes (f_1 f_2) + (g_1 , g_2) {f_2 d(f_1)} , 
\end{equation}
where $( \cdot , \cdot )$ is a symmetric non-degenerate invariant form on $\dg$ normalized by $(\theta ,\theta) = 2$,
where $\theta$ is the longest root of $\dg$.

 In agreement with the affine case, we add certain outer derivations to this algebra. More precisely, we consider
the following algebra of derivations:

\begin{equation}\label{D*}
\D^* = \sum_{p=1}^N \R d_p ,
\end{equation}
where $d_j = t_j\frac{\dd} {\dd t_j}$.
The derivations of $\R$ naturally extend to derivations of the Lie algebra $\dg\otimes \R$. Since the above algebra,
$\dg \otimes \R \oplus \K$, is the universal central extension of $\dg\otimes \R$ we can lift these derivations to this universal central extension by  using a result of Benkart-Moody
\cite{BM}, which says that every derivation of a perfect Lie algebra uniquely extends to a derivation of its universal central
extension. The action of $\D^*$ on $\dg \otimes \R \oplus \K$ is given by 
\begin{align}
[t^\m d_j, g\otimes t^\rb] &= r_j g\otimes t^{\m+\rb} ,\label {action1}\\
[t^\m d_j, t^\rb k_i] &= r_j t^{\m+\rb}k_i + \delta_{ji} \sum_{p=0}^N m_p t^{\m+\rb}k_p .\label {action2}
\end{align}
The first formula corresponds to the action of vector fields on functions, while the second is
the Lie derivative action of vector fields on 1-forms.

 The formulas \eqref{action1}, \eqref{action2} identify the bracket in $\D^*$ only up to a $\K$-valued 2-cocycle
$\tau \in H^2(\D^*, \K)$:
\begin{equation}\label{brackD}
[t^\m d_i , t^\rb d_j] = r_i t^{\m+\rb} d_j - m_j t^{\m+\rb} d_i + \tau( t^\m d_i , t^\rb d_j).
\end{equation}
We will use a  two-dimensional space of cocycles as in \cite{BB} (see the remarks between (2.12) and (2.13) there). It is spanned by the two cocycles
\begin{equation}\label{tau1}
\tau_1( t^\m d_i , t^\rb d_j) = m_j r_i \sum_{p=0}^N r_p t^{\m+\rb} k_p
= - m_j r_i \sum_{p=0}^N m_p t^{\m+\rb} k_p,
\end{equation}
and
\begin{equation}\label{tau2}
\tau_2( t^\m d_i , t^\rb d_j) = m_i r_j \sum_{p=0}^N r_p t^{\m+\rb} k_p
= - m_i r_j \sum_{p=0}^N m_p t^{\m+\rb} k_p.
\end{equation}
 Thus we obtain a two-parametric family of algebras $\g = \g_\tau = \dg \otimes \R \oplus \K \oplus \D^*$,
where the cocylce in \eqref{brackD} is a linear combination of \eqref{tau1}
and \eqref{tau2}:  $\tau = \mu \tau_1 + \nu \tau_2$. 
A consequence of the fact that the action \eqref{action2} of $\D^*$ on $\K$ is non-trivial, is that the center
of the algebra $\dg \otimes \R \oplus \K \oplus \D^*$ is finite-dimensional and is spanned
by $k_0, k_1, \ldots, k_N$.

As in the affine case, we can also consider a semidirect product of $\g$ with
the Virasoro algebra 
\begin{equation}
Vir = \C [t_0, t_0^{-1}] d_0 \oplus \C C_{Vir}.
\end{equation}
The Lie bracket in the Virasoro algebra is given by the well-known formula
\begin{equation}
[t_0^n d_0 , t_0^m d_0] = (m-n) t_0^{n+m} d_0 + \frac{n^3-n}{12} \delta_{n,-m}
C_{Vir} 
\end{equation}
and the element $C_{Vir}$ is central.

The action of $Vir$ on $\dg \otimes \R \oplus \K $ is given by \eqref{action1},
\eqref{action2}, where we let $j=0$ and $t^\m = t_0^m$, while
the action on $\D^*$ is similar to \eqref{brackD}, only we use 
a different cocycle:
\begin{equation}
[t_0^n d_0, t^\m d_i] = m_0 t_0 ^n t^\m d_i - \rho m_i n (n + 1) t_0 ^n t^\m k_0 , \quad \rho \in \C.
\end{equation}
This new cocycle emerged from the representations of $\g \oplus Vir$ that
we study below. Later we will see that we need to specify an exact value of 
$\rho$ and this will depend on the particular value of the central charge we 
are using.  

We next introduce the following formal fields in $(\g \oplus Vir)[[z,z^{-1}]]$ whose moments together with $C_{Vir}$ span $\g \oplus Vir$. Notice that in the single field
$K_0(\m,z)$ below the coefficient of $z^{-j-1}$ involves $t^{j+1}$ while in all other fields it only involves $t^j$. It turns out that this is important when we begin to line up our fields with elements in a VOA as it is tied up with the conformal weights, (see our comment below).
\begin{align}
g(\m,z)=&\sum _j g \otimes t_0^j t^\m  z^{-j-1}, \, g\in {\frak g}    \label{gm}\\
K_0(\m,z)=&\sum _j t_0^j t^\m k_0 z^{-j}  \label{Kzero}\\
K_s(\m,z)=&\sum _j t_0^j t^\m k_s z^{-j-1}  , \quad 1\leq s\leq N \label{Ka}\\
D_s(\m,z)=&\sum _j t_0^j t^\m d_s z^{-j-1}  , \quad 1\leq s\leq N  \label{Da}\\
L(z)=&-\sum _j t_0^j d_0 z^{-j-2}   \label{D0}
\end{align}
Here and in what follows, $\m \in \Z^N$ and $t^\m = t_1^{m_1}\ldots t_N^{m_N}$.

The Lie algebra structure of $\g \oplus Vir$ is encoded in the following 
relations for the above fields:
\begin{align}
\frac{d}{dz}K_0({\bf m},z)=&\sum_{p=1}^N m_p K_p({\bf m},z)\\
[g_1(\m,z_1),g_2(\rb,z_2)]=&\{[g_1,g_2](\m+\rb,z_2)+
(g_1,g_2)\sum_{p=1}^N m_pK_p(\m+\rb,z_2)\}z_1^{-1}\delta (\frac{z_2}{z_1})+\nonumber \\
&(g_1,g_2)K_0(\m+\rb,z_2)z_1^{-2}\delta^{(1)} (\frac{z_2}{z_1}) \label{g1g2}\\
[g(\m,z_1), K_i(\rb,z_2)] =& \quad 0 \label{gk}\\
[K_i(\m, z_1), K_j(\rb, z_2)] =& \quad 0 \label{kk}\\
[D_s(\m,z_1),g(\rb,z_2)]=& \quad r_sg(\m+\rb,z_2)z_1^{-1}\delta (\frac{z_2}{z_1})
\label{dagr}\\
[L(z_1),g(\m,z_2)]=&\left(\frac{\partial}{\partial z_2}(g(\m,z_2))\right)z_1^{-1}\delta (\frac{z_2}{z_1})+
g(\m,z_2)z_1^{-2}\delta^{(1)} (\frac{z_2}{z_1}) \label{d0gm}\\
[D_s(\m,z_1),K_0(\rb,z_2)]=&\quad r_sK_0(\m+\rb,z_2)z^{-1}_1\delta (\frac{z_2}{z_1})\label{dak0}\\
[D_s(\m,z_1),K_q(\rb,z_2)]=&
\{r_sK_q(\m+\rb,z_2)+\delta_{s,q}\sum_{p=1}^N m_p K_p(\m+\rb,z_2)\}z_1^{-1}\delta (\frac{z_2}{z_1}) \nonumber \\
&+ \delta_{s,q}K_0(\m+\rb,z_2)z_1^{-2}\delta^{(1)} (\frac{z_2}{z_1}) \label{dakb}\\
[L(z_1),K_0(\m,z_2)]=&\left(\frac{\partial}{\partial z_2} (K_0(\m,z_2))\right)z_1^{-1}
\delta (\frac{z_2}{z_1})\label{d0k0}\\
[L(z_1),K_s(\m,z_2)]=&\left(\frac{\partial}{\partial z_2}(K_s(\m,z_2))\right)z_1^{-1}\delta (\frac{z_2}{z_1})+
K_s(\m,z_2)z_1^{-2}\delta^{(1)} (\frac{z_2}{z_1}) \label{d0ka}\\
[D_s(\m,z_1),D_q(\rb,z_2)]=&
\{r_sD_q(\m+\rb,z_2)-m_qD_s(\m+\rb,z_2)\}z^{-1}_1\delta (\frac{z_2}{z_1})\nonumber\\
& - (\mu r_s m_q + \nu r_q m_s) \sum_{p=1}^Nm_p
K_p(\m+\rb,z_2)z^{-1}_1\delta (\frac{z_2}{z_1}) \nonumber \\
& - (\mu r_s m_q + \nu r_q m_s) K_0(\m+\rb,z_2)z_1^{-2} \delta^{(1)} (\frac{z_2}{z_1}) \label{dadb}
\end{align}
\begin{align}
[L(z_1),D_s(\m,z_2)]=&\frac{\partial}{\partial z_2}(D_s(\m,z_2))z_1^{-1}\delta (\frac{z_2}{z_1})+
D_s(\m,z_2)z_1^{-2}\delta^{(1)} (\frac{z_2}{z_1}) \nonumber \\
& + \rho  m_s K_0(\m,z_2) z_1^{-3}\delta^{(2)} (\frac{z_2}{z_1})
\label{d0da}\\
[L(z_1),L(z_2)]=&\left(\frac{\partial}{\partial z_2}(L(z_2))\right)z_1^{-1}\delta (\frac{z_2}{z_1})+2 L(z_2)
z_1^{-2}\delta^{(1)}(\frac{z_2}{z_1})+\nonumber\\
&\frac{C_{Vir}}{12}z_1^{-4}\delta ^{(3)}(\frac{z_2}{z_1}) \label{d0d0}
\end{align}
where $1\leq s,q\leq N, 0\leq i,j \leq N, \m,\rb\in\Z^N$ and $\delta ^{(n)}(z)$ is the n-th derivative of the delta series 
$\delta(z) = \sum_{j\in\Z} z^j$.  

Next we list the conformal weights of 
all the fields involved.  We recall that the conformal weights are defined relative to the operator $L(0)$ which appears in the Virasoro field:
\begin{equation} 
L(z)=\sum_{n\in\Z} L(n)z^{-n -2} \nonumber
\end{equation}
by declaring that the local field $b(z)$ has conformal weight $\Delta$ if 
\begin{equation}\label{conweight}
[L(0),b(z)]=(z\frac{d}{dz}+\Delta)b(z).  
\end{equation}
Since the field $L(z)$ satisfies the Virasoro commutation relations we 
will use the coefficient at the $-2 $nd power of $z$ in the expansion of 
$L(z)$.  Thus to read off the conformal weights so defined from \eqref{d0k0},
 \eqref{d0ka}, we extract the coefficient of $z_1^{-2}$.  We obtain that 
$K_0(\m,z)$ has conformal weight $0$, $D_s(\m,z),K_s(\m,z),g(\m,z)$ have conformal weight 
$1$ while $L(z)$ has conformal weight $2$.  

In order to prove that a map $\phi: \g \rightarrow {\text {End}}(V)$ is
a representation of $\g$, we would need to verify that the relations 
\eqref{g1g2}-\eqref{d0d0} hold in ${\text {End}}(V)$. 
We note that \eqref{dak0}-\eqref{d0ka} follow from \eqref{g1g2}-\eqref{d0gm}
since a derivation of a Lie algebra is uniquely determined by its action
on the generators, and the space $\dg\otimes \R$ generates the algebra
$\dg\otimes \R \oplus \K$. To make this more rigorous, we prove the following
simple Lemma:
\begin{lemma}\label{extder}
Let $\phi: G \rightarrow {\text End}(V)$ be a representation of a Lie algebra
$G$ and $d$ be a derivation of $G$. Let the subspace $X \subset G$ generate
$G$ as a Lie algebra. If $d^\prime \in {\text End} V$ satisfies
$\phi \circ d |_X = {\text ad}(d^\prime) \circ \phi |_X$ then

(1) $\phi \circ d = {\text ad}(d^\prime) \circ \phi $ on $G$ and

(2) The representation $\phi$ can be extended to a representation of a 
semi-direct product $\phi: G\oplus \C d \rightarrow  {\text End}(V)$ by
$\phi(d) = d^\prime$.
\end{lemma}

{\sl Proof.} Since $X$ generates $G$ as a Lie algebra, every element $g\in G$
may be written as $g = \sum_s [u_s, v_s] + w$, where $u_s, v_s, w \in X$. Then
\begin{equation}
\phi(d(g)) = \phi \left( \sum_s [d(u_s),v_s] + [u_s, d(v_s)] + d(w) \right) =
\nonumber
\end{equation}
\begin{equation}
\sum_s [\phi(d(u_s)),\phi(v_s)] + [\phi(u_s),\phi(d(v_s))] + \phi(d(w)) = 
\nonumber
\end{equation}
\begin{equation}
\sum_s [{\text {ad}}(d^\prime) \phi(u_s),\phi(v_s)] + [\phi(u_s),{\text ad}(d^\prime)\phi(v_s)] + {\text ad}(d^\prime)\phi(w) =   {\text ad}(d^\prime)\phi(g).
\nonumber
\end{equation}
This proves (1), while (2) immediately follows from (1).

\section{Vertex operator algebras}
We begin by reviewing some basic facts  about vertex operator algebras (VOAs) and begin by recalling the definition. For notation, we follow \cite{FHL}, \cite{FLM}.  
\begin{definition}\label{defVOA}
A vertex operator algebra is a 4-tuple (V,Y,${\bf 1}$,$\omega$) where $V=\oplus_{n \in \Z}V_n$ is a $\mathbb Z$-graded vector space 
equipped with a linear map 
\begin{equation*}
V\rightarrow (End\, V)[[z,z^{-1}]]: 
\end{equation*}
\begin{equation*}
v\rightarrow Y(v,z)=\sum_{n\in \mathbb Z}v_n z^{-n-1}, \qquad v_n \in 
End\quad V
\end{equation*}
and two distinguished vectors ${\bf 1}$ and $\omega \in V$ satisfying the following conditions:
\begin{align}
&u_nv=0,\,   \text{for $n>>0$  and $u,v\in V$ arbitrary},\\
&Y({\bf 1},z)=Id \\
&Y(v,z){\bf 1} \in V[[z]] \text{ and } v_{-1}{\bf 1}=v\\
&z_0^{-1}\delta(\frac{z_1-z_2}{z_0})Y(v,z_1)Y(w,z_2)-z_0^{-1}\delta(
\frac{z_2-z_1}{-z_0})Y(w,z_2)Y(v,z_1)\nonumber\\
&=z_2^{-1}\delta(\frac{z_1-z_0}{z_2})Y(Y(v,z_0)w,z_2) \quad \text (Jacobi \quad  identity), \\
&V_n = \{ 0\} \text{ for $n<<0$}
\end{align}
where as usual $\delta(z)=\sum_{n \in \Z} z^n$ and where  all binomial expressions are to be expanded in nonnegative powers of the second variable. The following Virasoro algebra relations must also hold:
\begin{equation}
[L(m),L(n)]=(m-n)L(m+n)+\frac{1}{12}(m^3-m)\delta _{m+n,0}(rank V) Id, 
\end{equation}
where  $m,n\in \mathbb Z$, $rank V \in \C$, $L(n)=\omega _{n+1}$ so that 
\begin{equation*}
Y(\omega,z)=\sum_{n \in \Z}\omega_n z^{-n-1}=\sum_{n \in \Z}L(n)z^{-n-2}
\end{equation*}
and
\begin{align}
L(0)v&=nv, \text{ for $v\in V_n$}\\
\frac{d}{dz}Y(v,z)&=Y(L(-1)v,z).
\end{align}
\end{definition}
This completes the definition.
By the above mentioned usual conventions we have that the three-variable delta function is given by:
\begin{equation}
\delta(\frac{z_1-z_2}{z_0})=\sum_{n\in \mathbb Z}z_0^{-n}(z_1-z_2)^n=\sum_{n\in 
\mathbb Z}\sum_{k\in \mathbb N}(-1)^k \binom{n}{k} z_0^{-n}z_1^{n-k}z_2^k
\end{equation}

Note that the above is not quite the definition of VOA as found in \cite {FHL} or \cite{FLM} as we have not required that the spaces $V_n$ be finite dimensional. This causes no problems in what follows and it is well known that most of the basic results in the theory of VOAs remain true without this extra hypothesis.

 The following important commutator formula is 
a consequence of the Jacobi identity:
\begin{equation}\label{comm}
\left[ Y(u,z_1), Y(v,z_2) \right] =
\sum_{n\geq 0} \frac{1}{n!} \left( z_1^{-1} \left( \frac{\dd}{\dd z_2} \right)^n
\delta \left( \frac {z_2}{z_1} \right) \right) Y(u_n v, z_2) ,
\end{equation}

Now we review the construction of the tensor product of VOAs that will be useful to us later.  
\eject
\begin{definition}\label{defTVOA}
Given finitely many vertex operator algebras 
$$
(V_1,Y_1,{\bf 1}_1,\omega_1),...,
(V_n, Y_n,{\bf 1}_n,\omega_n) 
$$ 
we call 
\begin{equation}
V=V_1\otimes \cdots \otimes V_n
\end{equation}
the tensor product of vertex algebras $V_1,...,V_n$ provided
\begin{align}
Y(v_1\otimes \cdots \otimes v_n,z)&=Y_1(v_1,z)\otimes \cdots \otimes Y_n(v_n,z), \quad v_i\in 
V_i\\
{\bf 1}&={\bf 1}_1\otimes \cdots {\bf 1}_n \\
\omega=\omega_1\otimes {\bf 1}_2\otimes &\cdots {\bf 1}_n+\cdots+ {\bf 1}_1\otimes 
\cdots \otimes {\bf 1}_{n-1}\otimes \omega _n
\end{align}
\end{definition}
Then the following theorem holds \cite{FHL}
\begin{theorem}\label{TVOA}
The tensor product $V$ of the vertex operator algebras 
$$
(V_1,Y_1,{\bf 1}_1,\omega_1,),...,
(V_n,Y_n,{\bf 1}_n,\omega_n,) 
$$
is again a vertex operator algebra. The rank of $V$ is equal 
to the sum of ranks of the $V_i$'s.  Moreover $V$ is a tensor product of Virasoro 
algebra modules.  
\end{theorem}

We next go on to review some particular VOA's. We begin with a lattice VOA and 
its subVOA which we need.

\noindent\begin{center}{\sl Lattice vertex algebras}
\end{center}
Let $\Lor$ be a free abelian group on $2N$ generators $\{a_i,b_i: 1\le i,j\le N\}$
We give $\Lor$ the structure of a {\sl lattice} by defining a symmetric bilinear
form:
\begin{align}
(.,.):\quad \Lor\times \Lor\rightarrow &\mathbb Z \text{ by }\nonumber\\
(a_i,b_j)=(a_i,b_j)=0& \qquad (a_i,b_j)=\delta _{i,j}
\end{align}
Note that the form $(.,.)$ is non-degenerate and $\Lor$  is even, 
that is $(x,x)\in 2\mathbb Z$.  
The construction of the VOA associated to $\Lor$ proceeds as follows.

First we consider the complexification of $\Lor$:
\begin{equation*}
H=\Lor \otimes_{\mathbb Z}{\mathbb C},
\end{equation*}
and extend $(.,.)$ by linearity obtaining a non-degenerate bilinear form on $H$.

We next ``affinize''  $H$ as follows. Thus, we define a Lie algebra: 
\begin{equation*}
\hat H=H\otimes {\mathbb C}[t,t^{-1}]\oplus {\mathbb C}K
\end{equation*}
with bracket:
\begin{equation*}
[x\otimes t^i,y\otimes t^j]=i(x,y)\delta_{i+j,0}K, \qquad [x\otimes t^i,K]=0.
\end{equation*}

We denote by $S_H$ the symmetric algebra of the space $H\otimes
t^{-1}{\mathbb C}[t^{-1}]$. Thus writing $\alpha(n)$ for the element $\alpha \otimes t^n$ for $\alpha \in H,n \in \Z$ we have that
\begin{equation*}  S_H=<\alpha_1(-n_1)..
.\alpha_k(-n_k)| \alpha_i \in H, n_i\in {\mathbb N}, 
  1\leq i \leq k>.
\end{equation*}
 Here, and throughout, we use pointed brackets, $<\quad >$ to denote the space spanned by the displayed vectors inside the brackets.

 We also need a twisted group algebra of $\Lor$, denoted ${\mathbb C}[\Lor]$, which we now describe. For this we use a two cocycle 
\begin{equation}\epsilon:\Lor \times \Lor \rightarrow \{ \pm1 \}
\end{equation}
which satisfies for all $\alpha,\beta,\gamma \in \Lor$
\begin{equation}
\epsilon(\alpha,\beta)\epsilon(\alpha+\beta,\gamma)=\epsilon(\beta,\gamma)\epsilon (\alpha, \beta +\gamma),
\end{equation}
\begin{equation}
\epsilon(\alpha,0)=\epsilon(0,\alpha)=1,
\end{equation}
and
\begin{equation}
\epsilon(\alpha,\beta)\epsilon(\beta,\alpha)=(-1)^{(\alpha,\beta)}.
\end{equation}

 The twisted group algebra of $\Lor$ which we need will be denoted $\C[\Lor]$ and it is spanned by elements $e^x$ with multiplication given by
\begin{equation}
e^\alpha e^\beta =\epsilon(\alpha,\beta)e^{\alpha+\beta}.
\end{equation}
We have $\C[\Lor]=<e^{x},x\in \Lor>$.  Then the space
$V_{\Lor}=S_H\otimes {\mathbb C}[\Lor]$ becomes the  VOA attached to the lattice $\Lor$.  The details can be found 
in \cite{FLM}. For example, we have
\begin{equation}
[a_p(n),b_q(m)]=n\delta_{p,q} \delta_{n+m,0}K,
\end{equation}
where $K$ acts as $1$.
  We will need an explicit form of $\omega$ as well as 
$L(-1)$.  With this in mind we briefly show below how these two elements are constructed.  
Let $\{\alpha_1,...,\alpha_{2N}\}$ be an orthonormal basis of $H$.   We set 
\begin{align}
{\bf 1}=&1\otimes 1,\\
\alpha(z)=&\sum _{n\in{\mathbb Z}}\alpha(n)z^{-n-1}\\
\omega=&\frac{1}{2}\sum _{k=1}^{2N} \alpha_k(-1)^2 \\
Y(\omega,z)=&\frac{1}{2}\sum _{k=1}^{2N} :\alpha_k(z)\alpha_k(z):\\
Y(1\otimes e^{\alpha},z)=&\exp(\sum _{n\ge 1}\frac{\alpha(-n)z^n}{n})
\exp(-\sum _{n\ge 1}\frac{\alpha(n)z^{-n}}{n})e^{\alpha} z^{\alpha}\label{YLe}
\end{align}
where $:\,:$ denotes  the normally ordered product and $e^{\alpha}$ acts by our twisted 
multiplication so that $e^{\alpha}e^{\beta}=\epsilon(\alpha,\beta)e^{\alpha + \beta}$.  Also $z^{\alpha}e^{\beta}=
z^{(\alpha,\beta)}e^{\beta}$.  
In particular from (3.22) one obtains, by expanding in powers of $z$, that 
\begin{equation}
L(m)=\frac{1}{2}\sum _{k=1}^{2N} \sum _{n\in {\mathbb Z}}:\alpha_k(m-n)\alpha_k(n):\,   \quad m\in {\mathbb Z}
\end{equation}

The rank of $V_{\Lor}$ is equal to $2N$.

One can give an explicit form of the map $Y$.  Set 
\begin{align}
\partial _n=&\frac{1}{n!}(\frac{d}{dz})^n\\
v=&\beta _1(-n_1)...\beta_k(-n_k) \otimes e^{\alpha}
\end{align}
Then 
\begin{equation}\label{YL}
Y(v,z)=:(\partial_{n_1-1}\beta _1(z)... \partial _{n_k -1}\beta _k(z))
Y(1\otimes e^{\alpha}):
\end{equation}
Consider the following {\sl sub-lattice}  $\Lor^+$ of $\Lor$ where $\Lor^+$ is generated by $<a_i,i=1,...,N>$.  
 We now define a subspace of $V_{\Lor}$ by saying
\begin{equation}\label{L^pm}
  V_{\Lor}^{+}=S_H\otimes {\mathbb C}[\Lor^{+}].
\end{equation}
We are going to need the fact that $V_{\Lor}^+$ is a sub-VOA of $V_{\Lor}$ and so we now indicate why this is true. 

The fact that $V_{\Lor}^+$ is a sub-VOA of $V_{\Lor}$ will follow from us showing that if $v \in V_{\Lor}^+$ then the moments of the series
$Y(v,z)$ map $V_{\Lor}^+$ to itself. Here the moments are just the $v_n$'s where $Y(v,z)=\sum_{n \in \Z}v_nz^{-n-1}$.  
However there are only three types of 
operators which act 
nontrivially on $1\otimes{\mathbb C}[\Lor]$:
$\{\alpha(0), \alpha \in \Lor\}$, $\{z^{\alpha},\alpha \in \Lor\} $ and the multiplication operators $\{e^{\alpha},\alpha \in \Lor\}$. From among them 
$\alpha(0)$ acts diagonally on
$1\otimes{\mathbb C}[\Lor]$, while $z^{\alpha}$ acts trivially on 
$1\otimes{\mathbb C}[\Lor^+]$ if $\alpha \in \Lor^+$. Finally, since 
$\Lor^+$ is a sub-lattice the 
twisted multiplication operator $e^{\alpha}$ maps $1\otimes{\mathbb C}
[\Lor^+]$ into itself if $\alpha \in \Lor^+$.  Finally, since \eqref{YL} contains
products of these three types of operators and operators which act trivially
on $1\otimes{\mathbb C}[\Lor]$ the fact that $V_{\Lor}^+$ is a sub-VOA of 
$V_{\Lor}$ follows.

\medskip 
One can also give an explicit form of the map $Y$ of this sub-VOA $V_{\Lor}^+ .$ Note that here, as usual, we are using the same symbol, $Y$, for both VOA's $V_{\Lor}$ and $V_{\Lor}^+$ but this should cause no confusion. In fact, for $v \in V_{\Lor}^+$ we see that
\begin{equation}
Y(v,z)=\sum_{n \in \Z}v_n|_{V_{\Lor}^+}z^{-n-1}.
\end{equation}
Since it is cumbersome to keep the notation for restriction we drop it and simply write
\begin{equation}
  Y(v,z)=\sum_{n \in \Z}v_n z^{-n-1}.
\end{equation}
Then using the same notation as in \eqref {YL} and assuming $\alpha \in \Lor^+$ we 
get:  
\begin{align}
Y(1\otimes e^{\alpha},z)=&\exp(\sum _{n\ge 1}\frac{\alpha(-n)z^n}{n})
\exp(-\sum _{n\ge 1}\frac{\alpha(n)z^{-n}}{n})e^{\alpha} \label{YLe+}\\
Y(v,z)=&:(\partial_{n_1-1}\beta _1(z)... \partial _{n_k -1}\beta _k(z))Y (1\otimes e^{\alpha},z): \label{YL+}
\end{align}
where now $e^{\alpha}$ acts by a standard multiplication.  Note also 
that since $z^{\alpha}$ acts trivially on $V_{\Lor}^+$, it disappeared from the 
formula for $Y(1\otimes e^{\alpha},z)$ (cf. \eqref{YLe}). 

  We also need the following notation as we will have need to deal with the vertex operators $Y(1\otimes e^\alpha,z)$ with $\alpha = \sum_{p=1}^N m_p a_p $.
 Thus for $\m=(m_1,...,m_N) \in \Z^N$ we let $\m \cdot {\bf a}=\sum_{p=1}^N 
m_p a_p$. That is, we think of ${\bf a}$ as the $N$-tuple $(a_1, ...,a_N).$ Then we have that
\begin{equation}
Y(1 \otimes e^{\sum_{p=1}^N m_p a_p},z)=Y(1 \otimes e^{\m \cdot {\bf a}},z).
\end{equation}
Also, in keeping with this notation we also let
\begin{equation}
\ma (-1)=\sum_{p=1}^N m_p a_p (-1).
 \end{equation}

\medskip

\begin{center} {\sl Affine vertex algebras}
\end{center}
For an affine VOA one lets $\gd$ be any finite dimensional simple Lie algebra and we let
$ (\cdot ,\cdot)$ be an invariant symmetric non-degenerate bilinear form on $\gd$ which is normalized so that 
$(\theta,\theta)=2$ where $\theta$ is the longest root of $\gd$. We let
 \begin{equation}
\widehat {\dg}=\dg \otimes \C[t_0,t_0^{-1}] \oplus \C k_0
 \end{equation}
be the non-twisted affine Lie algebra associated to this data with multiplication given by saying $k_0$ is central and
\begin{equation}
[x\otimes t_0^m,y \otimes t_0^n]=[x,y]\otimes t_0^{m+n}+m\delta_{m+n,0}(x,y)k_0.
\end{equation}
We write $x(n)$ for the element $x \otimes t_0^n$. Let $c$ be an arbitrary constant in $\C$ and let $\C 1$ be the one dimensional module for $\gd \otimes \C[t_0] \oplus \C k_0$ where $\gd \otimes \C[t_0]$ acts as zero on $\C 1$ and where $k_0$ acts as the constant $c$ on $1.$ We form a generalized Verma module by letting  
\begin{equation}
V_c = V_{\wdg, c} = {\text {Ind}}_{\gd \otimes \C[t_0] \oplus \C k_0}^{\wdg}\C 1
\end{equation}
This module has a unique irreducible quotient module which we denote by $L_c = L_{\wdg, c} .$ Then it is well known \cite {FZ}, \cite {Li} 
that as long as $c \neq -h^{\vee}$ then both $V_c $ and $L_c$ are VOA's. Here $-h^{\vee}$ is the negative of the dual Coxeter number of the Lie algebra $\gd.$

   We quickly review this so as to establish some notation. We work with $V_c$ first and define $Y(x(-1),z)=\sum_{n \in \Z}x(n)z^{-n-1} \in End(V_c)[[z,z^{-1}]]$ and let $D \in Der(\widehat \gd)$ be given by $[D, x(m)]=-mx(m-1),[D,k_0]=0.$ Then the set
\begin{equation}
\{ Y(x(-1),z) | x \in \dg \}
\end{equation}
consists of mutually local fields satisfying $\frac{d}{dz} Y(x(-1),z)=[D,Y(x(-1),z)]$ and so it follows that $V_c$ has the structure of a vertex algebra. Letting $x_1, ...,x_d$ be an orthonormal basis of $\gd$ we let
\begin{equation}
\omega=\frac{1}{2(c+h^\vee)}\sum_{i=1}^d x_i(-1)^2.
\end{equation}
As usual we let $Y(\omega,z)=\sum_{n \in \Z}L(n)z^{-n-2}$ and let 
\begin{equation}
(V_c)_n=\{v \in V_c | L(0)v=nv \}.
\end{equation} 
Then we have $V_c=\oplus_{n\in \Z}(V_c)_n$ and also
\begin{equation}
(V_c)_n=< y_1(-n_1)\cdots y_k(-n_k)| \sum_{i=1}^k n_i=n, n_i \in\N,y_i \in \gd >
\end{equation}
Then as long as $c \neq -h^{\vee}$ we obtain the VOA $(V_c, Y, {\bf 1}, \omega).$ We denote its rank by $r_c$ and have that 
\begin{equation}
r_c=\frac{c {\text{\ dim}}(\gd)}{c+ h^\vee}.
\end{equation}
   Finally, one knows that the usual radical of the module $V_c$ is  an ideal of the VOA $V_c$ so we obtain the VOA $L_c$ as a simple quotient of $V_c.$

\begin{center}{\sl $\gl$ vertex algebras}
\end{center}

The Lie algebra $gl_N$ is reductive and decomposes into a direct sum
$gl_N = sl_N \oplus \C I$. Let $\psi_1$ denote the projection on the traceless
matrices and $\psi_2$ be the projection on the scalar matrices in this
decomposition. Accordingly, we define the affine algebra $\gl$ to be a direct sum
of the affine algebra $\slN = sl_N(\C[t_0, t_0^{-1}]) \oplus \C C_1$ and a (degenerate) Heisenberg algebra $\h = I \otimes \C[t_0, t_0^{-1}] \oplus \C C_2$.

 The Lie bracket in $\gl$ is given by
\begin{align}
[g_1 \otimes t_0^i, g_2 \otimes t_0^j] = &[g_1,g_2] \otimes t_0^{i+j} 
+ i \delta_{i,-j} tr(\psi_1(g_1)\psi_1(g_2)) C_1 +\nonumber\\
&i\delta_{i,-j} (\psi_2(g_1),\psi_2(g_2)) C_2,\label{gl}
\end{align}  
where for the last term we will use the normalization $(I,I)=1$.

However, if we fix the basis of $gl_N$ consisting of elementary matrices,
then the formula above will become
\begin{align} 
[E_{pq} \otimes t_0^i, E_{rs} \otimes t_0^j] =&  \delta_{qr}E_{ps} \otimes t_0^{i+j}
-  \delta_{ps}E_{rq}\otimes t_0^{i+j}+\nonumber\\
&i \delta_{i,-j} \left(\delta_{qr}\delta_{ps} C_1 +
 +  \delta_{pq}\delta_{rs} \left( \frac{C_2}{N^2} -\frac {C_1}{N} \right)
 \right)\label{glN}
\end{align}  

To construct a VOA corresponding to $\gl$, we will take a tensor product
of a VOA for affine Lie algebra $\slN$ and a Heisenberg VOA corresponding to
$\h$. We have just seen the construction of an affine VOA and hence we obtain
 a Verma VOA $V_{\slN,c_1}$ and a simple VOA $L_{\slN,c_1}$ for an arbitrary complex
number $c_1 \neq -N$ (the dual Coxeter number for $sl(N)$ is $N$).
 Let us review the construction of the Heisenberg VOA.

Let $c_2$ be an arbitrary constant in $\C$ and
consider a one-dimensional module $\C 1$ for the subalgebra 
$I\otimes \C[t_0] \oplus <C_2>$ where $I\otimes \C[t_0]$ acts on $1$ trivially
and $C_2$ acts as multiplication by $c_2$.

 If $c_2 \neq 0$ then the induced module 
\begin{equation}
{\text{Ind}}_{I\otimes \C[t_0] \oplus <C_2>}^{\h} \C 1
\end{equation}
has a structure of a VOA \cite{K,MN}, which we denote by $L_{\h,c_2}$. 
As a vector space, $L_{\h,c_2}$ is isomorphic
to the symmetric algebra $S[I\otimes t_0^{-1}\C[t_0^{-1}]]$.
To construct the state-field correspondence $Y$, we take
\begin{align}
Y(I(-1),z) =& \sum_{n\in\Z} I(n)z^{-n-1}, \\
Y(I(-j-1),z) =&  \frac{1}{j!}\left( \frac{\dd}{\dd z} \right)^{j} 
Y(I(-1),z), \\
Y(I(-n_1)\ldots  I(-n_k),z) =& :Y(I(-n_1),z) \ldots Y(I(-n_k),z): ,
\end{align}
where we again use the notation $I(n) = I\otimes t_0^n$.

%

The Virasoro element is given by $\omega = \frac{1}{2 c_2} I(-1)^2$.

 If $c_2 = 0$ then the induced module has a one-dimensional factor
$L_{\h,0} = \C 1$. 
The VOA $L_{\h,c_2}$ is simple and its rank is equal to 1 if $c_2 \neq 0$
and to $0$ if $c_2 = 0$.

Finally, for a pair of complex numbers $c_1, c_2, \quad c_1 \neq -N$,
we define two VOAs for $\gl$ as tensor products: 
\begin{align}
V_{\gl} =& V_{\gl,c_1,c_2} := V_{\slN,c_1} \otimes L_{\h,c_2} \\
L_{\gl} =& L_{\gl,c_1,c_2} := L_{\slN,c_1} \otimes L_{\h,c_2}
\end{align}
Of course, both  VOAs $V_{\gl}$ and $L_{\gl}$ are modules for the Lie algebra $\gl$.

\medskip

\section{Toroidal VOAs}
Let $c, c_1, c_2$ be complex numbers such that $c \neq 0, c \neq -h^\vee,
c_1 \neq -N$.
We consider two tensor products of VOAs,$V=V_{\wdg,c}\otimes V_{\Lor}^+\otimes 
V_{{\gl,c_1,c_2}}$ and  $L=L_{\wdg,c}\otimes V_{\Lor}^+\otimes 
L_{{\gl,c_1,c_2}}$ with the following assignment of vertex operators,
which is the same for $V$ and $L$. We will take $\m \in \Z^N, 1 \leq s \leq N.$ Then we let: 

\begin{align}
K_0(\m,z)\rightarrow &c Y(1\otimes e^{\ma}\otimes 1,z)\label{vk0} \\
D_s(\m,z)\rightarrow &Y(1\otimes b_s(-1)e^{\ma}\otimes 1+
1\otimes e^{\ma}\otimes(\sum_{p=1}^Nm_p E_{ps}(-1)),z) \label{vda}\\
K_s(\m,z)\rightarrow&c Y(1\otimes a_s(-1)e^{\ma}\otimes 1,z) \label{vka}\\
g(\m,z)\rightarrow& Y(g(-1)\otimes e^{\ma}\otimes 1,z) \label{vg}\\
L(z)\rightarrow &Y(\omega_1\otimes 1\otimes 1+1\otimes \omega_2 \otimes 1+
1\otimes 1\otimes \omega_3,z) \label{vd0}\\
C_{Vir} \rightarrow& {\text { rank}}(L) id, \label{rank1}\\
where {\text { rank}}(L) =&
\frac{c {\text{dim}}(\gd)}{c+ h^\vee} + 2N + \frac{c_1 (N^2-1)}{c_1 + N} + 
{\text {rank}}(L_{Hei,c_2}) \label{rank}
\end{align}
We will denote by $Y_i,i=1,2,3$ the restriction of the state-field map $Y$ 
to the first, second, and third factor respectively in the tensor product.  
\begin{theorem} \label{main} 
Under the correspondence \eqref{vk0}-\eqref{rank}, both $V$ and $L$ are modules for
the Lie algebra $\g_\tau \oplus Vir = \dg\otimes\R \oplus \K \oplus \D^* \oplus Vir$,
where $\tau = \mu\tau_1 + \nu \tau_2$ with $\mu = \frac{1 - c_1}{c}$, 
$\nu = 
\frac{\frac{c_1}{N} - \frac{c_2}{N^2}}{c}$ and $\rho =\frac{1}{2c}$.  Moreover,
 $C_{Vir}$ equals {\text { rank}}($L$) as given by \eqref{rank}. 

\end{theorem}
{\sl Proof.}

\noindent 
The proof is identical for $V$ and $L$, so we will carry it out for $V$.
Our first claim is that 
\begin{equation}\label{constraint}
\frac{d}{dz}Y(1\otimes e^{\ma}\otimes 1,z)=
\sum_{p=1}^N m_p Y(1\otimes
a_p (-1)e^{\ma}\otimes 1,z)
\end{equation}
It suffices to prove this claim for $Y_2$.  From the VOA axioms we have 
\begin{equation}
\frac{d}{dz}Y_2(e^{\ma},z)=Y_2(L(-1)e^{\ma},z). 
\end{equation} 
To compute 
$L(-1)e^{\ma}$ we observe that 
\begin{align}
L(-1)e^{\ma}=&\lim_{w\rightarrow 0}[L(-1),Y_2(e^{\ma} ,w)]1=\nonumber \\
&\lim_{w\rightarrow 0}\frac{d}{dw}Y_2(e^{\ma},w)1
\end{align}
>From \eqref {YLe} we obtain $\sum_{p=1}^N m_p a_p(-1)e^{\ma}$ for the 
limit and the claim follows.

In order to prove that $V$ is a module for the toroidal Lie algebra,
we need  to show that the relations \ref{g1g2} - \ref{d0d0} hold in
End$(V)[[z,z^{-1}]]$.
We begin with the commutator corresponding to \eqref{g1g2}.

\noindent The first step here is to verify all commutators involving the 
affine part.  We observe that we can ignore the the third 
factor in the tensor product.  So we will suppress the third factor 
and will restore it only in the final answers. First,
in our affine VOA we have the following relation:
\begin{align}
[Y_1(g_1(-1),z_1),Y_1(g_2(-1),z_2)] =&
Y_1([g_1,g_2](-1),z_2)z_1^{-1}\delta(\frac{z_2}{z_1}) + \nonumber \\
&z_1^{-1} (g_1,g_2)c\frac{\partial}{\partial z_2}\delta(\frac{z_2}{z_1})
\end{align}
 
Using this we get 
\begin{align}
&[Y_1(g_1(-1),z_1)\otimes Y_2(e^{\ma},z_1),Y_1(g_2(-1),z_2)\otimes Y_2(e^{\ra},z_2)]=\nonumber \\
&\{Y_1([g_1,g_2](-1),z_2)z_1^{-1}\delta(\frac{z_2}{z_1}) + 
z_1^{-1} (g_1,g_2)c\frac{\partial}{\partial z_2}\delta(\frac{z_2}{z_1})\}\otimes Y_2(e^{\ma},z_1)Y_2(e^{\ra},z_2)=\nonumber\\
&Y_1([g_1,g_2](-1),z_2)\otimes Y_2(e^{\mra },z_2)z_1^{-1}\delta(\frac{z_2}{z_1}) + \nonumber \\
&(g_1,g_2)z_1^{-1}c\frac{\partial}{\partial z_2}\delta(\frac{z_2}{z_1})\otimes Y_2(e^{\ma},z_1)Y_2(e^{\ra},z_2)
\end{align}
The second term above can be dealt with using properties of the 
$\delta$-function, yielding 
\begin{align}
&[Y_1(g_1(-1),z_1)\otimes Y_2(e^{\ma},z_1),Y_1(g_2(-1),z_2)\otimes Y_2(e^{\ra},z_2)]= \nonumber \\
&\{Y_1([g_1,g_2](-1),z_2)\otimes Y_2(e^{\mra },z_2)+(g_1,g_2)c
1\otimes Y_2(\ma (-1)e^{\mra },z_2)\}z_1^{-1}\delta(\frac{z_2}{z_1}) + \nonumber \\
&(g_1,g_2)z_1^{-1}c 1 \otimes Y_2(e^{\mra },z_2)\frac{\partial}{\partial z_2}\delta(\frac{z_2}{z_1})
\end{align}
which, in turn, implies 
\begin{align}\label {vg1g2}
&[Y(g_1(-1)\otimes e^{\ma}\otimes 1,z_1),
Y_2(g_2(-1)\otimes e^{\ra}\otimes 1,z_2)]=\nonumber \\
&\{Y([g_1,g_2](-1)\otimes e^{\mra }\otimes 1,z_2)+(g_1,g_2)c Y(1\otimes \ma (-1)e^{\mra }\otimes
1, z_2)\}z_1^{-1}\delta(\frac{z_2}{z_1}) + \nonumber \\
&(g_1,g_2)z_1^{-1}c Y(1\otimes e^{\mra }\otimes 1, z_2)\frac{\partial}{\partial z_2}\delta(\frac{z_2}{z_1})
\end{align}
thus reproducing \eqref{g1g2}.  

 To establish \ref{kk} we note that the operators
$Y_2(e^\ma,z)$ and $Y_2(a_s(-1) e^\ma, z)$ correspond to elements from the totally
isotropic subspace in $\Lor$ and thus all commute. This also implies that
\begin{equation}
[Y_1(g(-1),z_1)\otimes Y_2(e^\ma,z_1) \otimes 1,
1 \otimes Y_2(e^\ra,z_2) \otimes 1] = 0
\end{equation}
and
\begin{equation}
[Y_1(g(-1),z_1)\otimes Y_2(e^\ma,z_1) \otimes 1,
1 \otimes Y_2(a_s (-1) e^\ra,z_2) \otimes 1] = 0,
\end{equation}
which establishes \eqref{gk}.

We turn to verifying \eqref{dagr}.  We start with 
\begin{align*}
&[1\otimes Y_2(\bs (-1)e^{\ma},z_1)\otimes 1 +
1\otimes Y_2(e^{\ma},z_1)\otimes Y_3(\sum _{p=1}^N m_p E_{ps}(-1),z_1),\\
&Y_1(g(-1),z_2)\otimes Y_2(e^{\ra},z_2)\otimes 1]=\\
&Y_1(g(-1),z_2)\otimes [Y_2(\bs (-1)e^{\ma},z_1),Y_2(e^{\ra},z_2)]\otimes 1
\end{align*}

 We compute the commutator $[Y_2(e^{\ra},z_2), Y_2(\bs (-1)e^{\ma},z_1)]$
using the commutator formula \eqref{comm}:
\begin{align}
[Y_2(e^{\ra},z_2), Y_2(\bs (-1)e^{\ma},z_1)] = \nonumber\\
\sum_{n\geq 0} \frac{1}{n!} \left( z_2^{-1} \left( \frac {\dd}{\dd z_1}
\right)^n \delta \left( \frac {z_1}{z_2} \right) \right) 
Y\left( e^{\ra}{}_{(n)} b_s(-1) e^\ma , z_1 \right) .\label{com1}
\end{align}
The terms $e^{\ra}_{(n)} b_s(-1) e^\ma, \quad n\geq 0,$ are given by the 
antiholomorphic part of $Y_2(e^{\ra},z) b_s(-1) e^\ma$. Here is what we mean by this.
For an arbitrary expression $x(z) = \sum_{k\in\Z} x_k z^k \in {\text{End}}(V)[[z,z^{-1}]]$,
we denote its antiholomorphic part $\sum_{k<0} x_k z^k$ by $x(z)^-$.
Since
\begin{equation}
Y_2(e^{\ra},z) b_s(-1) e^\ma = [Y_2(e^{\ra},z), b_s(-1)] e^\ma
+ b_s(-1) Y_2(e^{\ra},z)e^\ma
\end{equation}
and the antiholomorphic part of the last term is 0,
while the commutator in the first term can be easily computed using 
the explicit expression \eqref{YLe} for $ Y_2(e^{\ra},z)$:
\begin{equation}\label{eb}
[Y_2(e^{\ra},z), b_s(-1)] = -\frac{r_s}{z} Y_2(e^{\ra},z),
\end{equation}
we get that
\begin{equation}
\left( Y_2(e^{\ra},z) b_s(-1) e^\ma \right)^- = \left(
 -\frac{r_s}{z} Y_2(e^{\ra},z) e^\ma \right)^- = 
-\frac{r_s}{z} e^\mra .
\end{equation}
This means that $ e^{\ra}_{(0)} b_s(-1) e^\ma = - r_s e^\mra $ and
 $ e^{\ra}_{(n)} b_s(-1) e^\ma = 0$ for $n > 0$. Substituting this into
\eqref{com1} we obtain the desired commutator relation \eqref{dagr}:
\begin{align}
&[Y(1\otimes b_s(-1)e^{\ma}\otimes 1+
1\otimes e^{\ma}\otimes(\sum_{p=1}^Nm_p E_{ps}(-1)),z_1),
Y(g(-1)\otimes e^{\ra}\otimes 1,z_2)] =\nonumber \\
& r_s z_2^{-1} \delta \left( \frac {z_1}{z_2} \right) 
Y_1(g(-1),z_2)\otimes Y_2(e^\mra, z_1) \otimes 1 = \nonumber \\
& r_s z_1^{-1} \delta \left( \frac {z_2}{z_1} \right) 
Y(g(-1)\otimes e^\mra \otimes 1, z_2).
\end{align}

We turn now to \eqref{d0gm}:
\begin{align}
& [Y(\omega,z_1), Y(g(-1)\otimes e^\ma\otimes 1, z_2)] = \nonumber\\
& \sum_{n\geq 0} \frac{1}{n!} \left( z_1^{-1} \left( \frac {\dd}{\dd z_2}
\right)^n \delta \left( \frac {z_2}{z_1} \right) \right) 
Y\left( \omega_{n} (g(-1) \otimes e^\ma \otimes 1) , z_2 \right) . \label{Lg}
\end{align}

Using the general formula for the weight of the n-th product in a VOA
\begin{equation}\label{wt}
wt(u_n v) = wt(u) + wt(v) - n - 1
\end{equation}
 and using the fact that
$wt(\omega) = 2$ and $wt( g(-1) \otimes e^\ma \otimes 1 ) = 1$ we see that
$\omega_n (g(-1) \otimes e^\ma \otimes 1) = 0$ for $n \geq 3$.
We thus need to compute $\omega_n (g(-1) e^\ma \otimes 1) = L(n-1) (g(-1)\otimes e^\ma \otimes1)$, for
$n = 0, 1, 2$.

 Observing that $L(0) v = nv$ for $v\in V_n$, and taking into account
that $wt( g(-1) \otimes e^\ma \otimes 1) = 1$, we get 
$L(0) (g(-1) \otimes e^\ma \otimes 1) = g(-1) \otimes e^\ma \otimes 1$.

 For $L(-1)$ we use the relation $Y(L(-1)v,z) = \frac {\dd}{\dd z} Y(v,z)$.

 To compute the term with $L(1)$ we consider
\begin{align}
L(1) (g(-1) \otimes e^\ma \otimes 1) =& 
(L(1) g(-1)) 1 \otimes  e^\ma \otimes 1 + g(-1) \otimes L(1) e^\ma \otimes 1 \nonumber\\
& g(-1) \otimes e^\ma \otimes L(1)1.
\end{align}
However all three terms here are zero. The last two terms vanish because of
\eqref{wt}, while for the first term we use the following relation in the affine VOA:
\begin{equation}
[L(n), g(k)] = -k g(n+k), \quad {\text for \ } g\in\dg ,
\end{equation}
which gives $[L(1), g(-1)] 1 = g(0) 1 = 0$.

Summarizing, we get the relation \eqref{d0gm}:
\begin{align}
[L(z_1), Y(g(-1)\otimes e^\ma\otimes 1, z_2)]=&
z_1^{-1}\delta (\frac{z_2}{z_1}) \left(\frac{\partial}{\partial z_2}
 Y(g(-1)\otimes e^\ma\otimes 1, z_2)\right)+\nonumber\\
&z_1^{-1} \left( \frac {\dd}{\dd z_2}
\right) \delta (\frac{z_2}{z_1}) Y(g(-1)\otimes e^\ma\otimes 1, z_2). 
\end{align}

 By Lemma \eqref{extder}, the relations \eqref{dak0} - \eqref{d0ka} follow
automatically from the ones we just established.

To verify \eqref{dadb} we write out the whole commutator: 
\begin{align*}
&[1\otimes Y_2(\bs(-1)e^{\ma},z_1)\otimes 
1 + 1\otimes Y_2(e^{\ma},z_1)\otimes Y_3(
\sum _{p=1}^N m_j E_{ps}(-1),z_1),\\
&1\otimes Y_2(\bq(-1)e^{\ra},z_2)\otimes 
1 + 1\otimes Y_2(e^{\ra},z_2)\otimes Y_3(
\sum _{k=1}^N r_k E_{kq}(-1),z_2)]=\\
&1\otimes [Y_2(\bs(-1)e^{\ma},z_1),Y_2(\bq(-1)e^{\ra},z_2)]\otimes 1+\\
&1\otimes [Y_2(\bs(-1)e^{\ma},z_1),Y_2(e^{\ra},z_2)]
\otimes Y_3(\sum _{k=1}^N r_k E_{kq}(-1),z_2)+\\
&1\otimes [Y_2(e^{\ma},z_1),Y_2(\bq(-1)e^{\ra },z_2)]\otimes Y_3(\sum _{p=1}^N m_p E_{ps}(-1),z_1)\\
&1\otimes Y_2(e^{\ma},z_1)Y_2(e^{\ra},z_2)\otimes 
[Y_3(\sum _{p=1}^N m_p E_{ps}(-1),z_1),Y_3(\sum _{k=1}^N r_k 
E_{kq}(-1),z_2)],
\end{align*}
where in the very last line we used that 
$[Y_2(e^{\ma},z_1),Y_2(e^{\ra},z_2)]=0$.  
To compute the first commutator above we will need 
$(Y_2(\bs(-1)e^{\ma},z)\bq (-1)e^{\ra})
^-$
We claim that 
\begin{align}
&(Y_2(\bs(-1)e^{\ma},z)\bq (-1)e^{\ra})^-= \nonumber\\
& \frac{1}{z}(r_s\bq(-1)-m_q \bs(-1) )e^{\mra}  
- \frac{1}{z} r_s m_q \sum_{p=1}^N m_p a_p (-1) e^\mra
- \frac{1}{z^2} r_s m_q e^\mra \label{bambbr}
\end{align}

\noindent To prove this statement we write: 
\begin{align}
&(Y_2(\bs(-1)e^{\ma},z)\bq (-1)e^{\ra})^-=
(:\bs (z)Y_2(e^{\ma},z): \bq(-1)e^{\ra})^-=\nonumber\\
&(\bs (z)^+ Y_2(e^{\ma},z)\bq(-1)e^{\ra}1)^-
+(Y_2(e^{\ma},z)\bs (z)^- \bs(-1)e^{\ra }1)^-
\label{first}
\end{align}
Now the first term in \eqref{first}
\begin{align}
&(\bs (z)^+ Y_2(e^{\ma},z)\bq(-1)e^{\ra}1)^-= \nonumber \\
&(\bs(z)^+[Y_2(e^{\ma},z),\bq(-1)]e^{\ra}1)^-+
(\bs(z)^+\bq(-1)Y_2(e^{\ma},z)e^{\ra} 1)^-
\end{align}
Since the second term is $0$ we get using \eqref{eb}:
\begin{align}
(\bs (z)^+ Y_2(e^{\ma},z)\bq(-1)e^{\ra}1)^-&=(\bs(z)^+[Y_2(e^{\ma},z),\bq(-1)]e^{\ra}1)^- \nonumber\\
&= - \left( \frac{m_q}{z}\bs(z)^+ Y_2(e^{\ma},z) e^\ra 1\right)^-,
\end{align}
which in turn, after noticing that only 
$\lim_{z\rightarrow 0}\bs(z)^+=\bs(-1)$ 
contributes to the antiholomorphic part, gives
\begin{equation}
(\bs (z)^+ Y_2(e^{\ma},z)\bq(-1)e^{\ra}1)^-=
-\frac{m_q \bs(-1)e^{\mra }}{z}
\end{equation}
Likewise, to compute the second term in \eqref{first}, we make use of 
the commutator
\begin{equation} \label{be}
[\bs(z)^-, e^{\ra}]=\frac{r_s e^{\ra}}{z}
\end{equation}
as well as of \eqref{eb}:
\begin{align}
&(Y_2(e^{\ma},z)\bs (z)^- \bq(-1)e^{\ra }1)^- =
(Y_2(e^{\ma},z) \bq(-1) \frac{r_s}{z}e^{\ra }1)^- = \nonumber\\
&(\frac{r_s}{z} [Y_2(e^{\ma},z) \bq(-1)] e^{\ra }1)^- +
(\frac{r_s}{z} \bq(-1) Y_2(e^{\ma},z)   e^{\ra }1)^- = \nonumber\\
&-( \frac{r_s m_q}{z^2} Y_2(e^{\ma},z)   e^{\ra }1)^-
+ (\frac{r_s}{z} \bq(-1) Y_2(e^{\ma},z)   e^{\ra }1)^- = \nonumber\\
& - \frac{r_s m_q}{z^2} e^\mra -  \frac{r_s m_q}{z} \sum_{p=1}^N m_p a_p(-1) e^\mra
+ \frac{r_s}{z}  \bq(-1) e^\mra ,
\end{align}
from which the  \eqref{bambbr} follows.

\noindent Using \eqref{be} we get 
\begin{equation} \label{secondo}
(Y_2(\bs(-1)e^{\ma},z)e^{\ra})^-=\frac{r_s e^{\mra}}{z}
\end{equation}
and, using \eqref{eb} we get
\begin{equation} \label {third}
(Y_2(e^{\ma},z)\bq(-1)e^{\ra})^-=-\frac{m_q e^{\mra}}{z}
\end{equation}

\noindent There is one more commutator, namely,
\begin{align*}
1\otimes Y_2(e^{\ma},z_1)Y_2(e^{\ra},z_2)
\otimes [Y_3(\sum_{p=1}^N m_p E_{ps}(-1),z_1),Y_3(\sum _{k=1}^N
r_kE_{kq}(-1),z_2)]
\end{align*} 
Using that $V_{\hat gl(N)}$ is a ${\hat gl(N)}$ module we obtain
\begin{align*}
&[Y_3(\sum_{p=1}^N m_p E_{ps}(-1),z_1),Y_3(\sum _{k=1}^N
r_kE_{kq}(-1),z_2)]=\\
&z_1^{-1} \delta (\frac{z_2}{z_1})\{r_sY_3(\sum_{p=1}^N
m_p E_{pq}(-1),z_2)-m_q Y_3(\sum_{k=1}^N r_k E_{ks}(-1),z_2)\}+\\
&(r_s m_q c_1 + r_q m_s \left( \frac{c_2}{N^2} - \frac{c_1}{N}\right)) z_1^{-1}\frac{\partial}{\partial z_2} \delta(\frac{z_2}{z_1})
\end{align*}
We compute the contribution of
the very last term, leaving out the factor
$r_s m_q c_1 + r_q m_s \left( \frac{c_2}{N^2} - \frac{c_1}{N}\right)$, 
with the help of the identity $
z_1^{-1}\frac{\partial}{\partial z_2} \delta(\frac{z_2}{z_1})=-
z_2^{-1}\frac{\partial}{\partial z_1} \delta(\frac{z_2}{z_1})$ 
\begin{align}\label{lastterm}
&1\otimes Y_2(e^{\ma},z_1)Y_2(e^{\ra},z_2) z_1^{-1}\frac{\partial}{\partial z_2} \delta(\frac{z_2}{z_1})\otimes 1= \nonumber\\
&-1\otimes Y_2(e^{\ma},z_1)Y_2(e^{\ra},z_2) z_2^{-1}\frac{\partial}{\partial z_1} \delta(\frac{z_2}{z_1})\otimes =\nonumber\\ 
-&1\otimes (\frac{\partial}{\partial z_1}(Y_2(e^{\ma},z_1)Y_2(e^{\ra},z_2) z_2^{-1}\delta(\frac{z_2}{z_1})))\otimes 1 +\nonumber\\
&z_2^{-1} 1\otimes Y_2(e^{\ra},z_2)(\frac{\partial}{\partial z_1}
Y_2(e^{\ma},z_1)) \delta(\frac{z_2}{z_1}))\otimes 1=\nonumber\\
&- 1 \otimes Y_2(e^{\mra },z_2)z_2^{-1}\frac{\partial}{\partial z_1}\delta(\frac{z_2}{z_1})\otimes 1+\nonumber\\
&z_2^{-1} 1\otimes Y_2(e^{\ra},z_2)
\sum_{p=1}^N m_p Y_2(a_p (-1) e^\ma, z_1)
\delta(\frac{z_2}{z_1}))\otimes 1=\nonumber\\
& 1 \otimes Y_2(e^{\mra },z_2)z_1^{-1}\frac{\partial}{\partial z_2}\delta(\frac{z_2}{z_1})\otimes 1+\nonumber\\
& z_2^{-1} 1\otimes 
\sum_{p=1}^N m_p Y_2(a_p (-1) e^\mra, z_2)
\delta(\frac{z_2}{z_1}))\otimes 1 
\end{align}
Here we used the fact that $Y_2(e^{\ra},z_1)
 Y_2(e^\ma, z_2) z_1^{-1} \delta(\frac{z_2}{z_1}) =
 Y_2(e^\mra, z_2) z_1^{-1} \delta(\frac{z_2}{z_1})$.

For simplicity, we set 
\begin{equation}\label{daz}
d_s(\m ,z)=Y(1\otimes \bs(-1)e^{\ma}\otimes 1+
1\otimes e^{\ma}\otimes(\sum_{p=1}^N m_p E_{ps}(-1)),z)
\end{equation}
Combining together \eqref {bambbr}, \eqref{secondo}, \eqref{third}, 
\eqref{lastterm} we get that 
\begin{align}\label {vdadb}
&[d_s(\m ,z_1),d_q(\rb ,z_2)]=\nonumber\\
&\{r_s d_q({\bf m+r},z_2)-m_q d_s({\bf m+r},z_2) + \nonumber\\
&\left(r_s m_q (c_1-1) + r_q m_s \left( \frac{c_1}{N} - \frac{c_2}{N^2}\right)
\right)
\sum_{p=1}^N m_p Y(1\otimes a_p(-1) e^{\mra }\otimes 1,z_2)\}
z_1^{-1}\delta(\frac{z_2}{z_1})+ \nonumber\\
&\left(r_s m_q (c_1-1) + r_q m_s \left( \frac{c_1}{N} - \frac{c_2}{N^2}\right)
\right)
Y(1\otimes e^{\mra }\otimes 1,z_2)z_1^{-1}\frac{\partial}{\partial z_2}\delta(\frac{z_2}{z_1})
\end{align}
which coincides with \eqref{dadb} with $\mu = \frac{1 - c_1}{c}$ and 
$\nu = \frac{\frac{c_1}{N} - \frac{c_2}{N^2}}{c}$ .

\medskip

To verify \eqref {d0da} we need to inspect the following 
commutators
\begin {align*}
&[Y_2(\omega _2,z_1),Y_2(\bs(-1)e^{\ma},z_2)]\\ 
&[Y_2(\omega _2,z_1),Y_2(e^{\ma},z_2)]\\
&[Y_3(\omega _3,z_1),Y_3(\sum _{p=1}^N m_pE_{ps}(-1),z_2)]
\end{align*}
In all three cases the verification 
is similar to \eqref{Lg}.  We use the fact that 
the weight of $\bs(-1)e^{\ma}$ and $E_{ps}(-1)$ is one, and is zero for 
$e^{\ma}$. 
The results for the anti-holomorphic parts have the same form in all three cases
\begin{align*}
&(Y(\omega,z)v)^-=\frac{L(-1)v}{z}+\frac{L(0)v}{z^2}+\frac{L(1)v}{z^3}
\end{align*}
and in a way analogous to the computation for \eqref{Lg}, we obtain
\begin{align*}
&(Y_2(\omega,z)\bs(-1)e^{\ma})^-=
\frac{1}{z}L(-1)\bs(-1)e^{\ma}+\frac{1}{z^2}\bs(-1)e^{\ma} + 
\frac{m_s}{z^3}e^{\ma}, \\
&(Y_2(\omega,z)e^{\ma})^-=
\frac{1}{z}L(-1)e^{\ma},\\
&(Y_3(\omega,z)\sum _{p=1}^Nm_p E_{ps}(-1))^-=
\frac{1}{z}L(-1)\sum _{p=1}^N m_pE_{ps}(-1)+\frac{1}{z^2}\sum _{p=1}^N m_pE_{ps}(-1)
\end{align*}
which after assembling all the terms yields \eqref {d0da}, provided 
$\rho=\frac{1}{2c}$.

\medskip 

Finally, the field $Y(\omega,z)$ satisfies the Virasoro commutation relations 
by definition, with the proviso $C_{Vir}={\text {rank}}(L)$.

\section{Modules for toroidal VOAs}

 In \cite{BB} a family of irreducible $\g$-modules was constructed.
In this section we will identify these modules as irreducible modules
for the toroidal VOAs, thus giving their explicit realizations.

We begin by recalling the construction of the $\g_\tau$ modules from \cite{BB}.
Fix as before a non-zero complex number $c$. The input data for building 
the module consists of a finite-dimensional irreducible $\dg$-module
$\dV$, a finite-dimensional irreducible $sl_N$-module $W$ and a complex 
number $c_3$. The pair $(W,c_3)$ actually defines a $gl_N$-module structure
on the space $W$ by letting the identity matrix act as multiplication by $c_3$.

The toroidal Lie algebra $\g = \dg\otimes\R \oplus \K \oplus \D^*$ has a 
$\Z$-grading by  degree in the variable $t_0$. We also consider a decomposition
$\g = \g_- \oplus \g_0 \oplus \g_+$ associated with this grading. The 
subalgebra $\g_0$ is spanned by the elements 
$\left\{ g \otimes t^\m , t^\m k_0, t^\m k_s, t^\m d_s | \m \in\Z^n, g\in\dg, 
s = 1,\ldots,N \right\}$. 
We begin by constructing a $\g_0$-module $T$ from the input data:
\begin{equation}
T =  \dV \otimes \C[\Lor^+]\otimes W.
\end{equation}
The action of $\g_0$ on $T$ is given by the formulas (see \cite{BB}):
\begin{align}
g \otimes t^\m  \left( v\otimes e^\ra \otimes w \right) = 
& (gv)\otimes e^\mra \otimes w, \label{act1}\\
t^\m k_0 \left( v\otimes e^\ra \otimes w \right) = 
&  c v\otimes e^\mra \otimes w, \\
t^\m k_s \left( v\otimes e^\ra \otimes w \right) = 
&  0, \\
t^\m d_s \left( v\otimes e^\ra \otimes w \right) = 
&  r_s v\otimes e^\mra \otimes w 
+ \sum_{p=1}^N m_p v\otimes e^\mra \otimes E_{ps} w , \label{act2}
\end{align}

Next we let $\g_+$ act trivially on $T$  and we build the induced module
\begin{equation}
\overline U = {\text {Ind}}_{\g_0 \oplus \g_+}^{\g} (T) .
\end{equation}
The module $\overline U$ inherits a natural $\Z^{N+1}$-grading from $\g$. 
We point out that the $\Z^{N+1}$-homogeneous components of $\overline U$
lying below the top $T$ are infinite-dimensional. The module $\overline U$
has a unique maximal homogeneous submodule  $\overline U_{\text{rad}}$ intersecting
$T$ trivially. The factor-module
\begin{equation}
U = U(c,\dV,W,c_3) = \overline U / \overline U_{\text {rad}}
\end{equation}
is a graded simple module (every $\Z^{N+1}$-graded submodule is trivial).
Applying Theorem 1.12 of \cite{BB}, we get that the homogeneous components 
of $U$ are finite-dimensional.

 Now we will use the same input data to construct a module for the toroidal
VOA $V_{c,c_1,c_2}$, where $c$ is as above and $c_1,c_2$ are related to
the cocycle $\tau$ by $\mu = \frac{1 - c_1}{c}$, 
$\nu = \frac{\frac{c_1}{N} - \frac{c_2}{N^2}}{c}$. We also impose additional restrictions, namely 
$c \neq -h^\vee, c_1 \neq -N$. The toroidal VOA $V_{c,c_1,c_2}$ is a
tensor product of three VOAs:  $V_{c,c_1,c_2} = V_{\wdg,c}\otimes V_{\Lor}^+
\otimes V_{\gl,c_1,c_2}$. We will construct a VOA module by taking a tensor
 product of modules for each factor:
\begin{equation}
M_{\wdg,c,\dV} \otimes V_{\Lor}^+ \otimes M_{\gl, c_1,c_2,c_3,W}
\end{equation}
 We view $V_{\Lor}^+$ as a module over itself. Let us briefly describe
$M_{\wdg,c,\dV}$ and $M_{\gl, c_1,c_2,c_3,W}$.

 Extend the action of $\dg$ on $\dV$ to an action of $\dg \otimes \C[t_0]
\oplus \C k_0$ by letting  $\dg \otimes t_0\C[t_0]$ act on $\dV$ trivially
and letting $k_0$ act as multiplication by $c$. Consider the induced
$\wdg$-module:
\begin{equation}
\overline M_{\wdg,c,\dV} =
{\text {Ind}}_{\dg \otimes \C[t_0]\oplus \C k_0}^{\wdg}
(\dV)
\end{equation}
and its irreducible factor $M_{\wdg,c,\dV}$. The $\wdg$-module
$M_{\wdg,c,\dV}$ is also an irreducible VOA module for affine VOA $V_{\wdg,c}$.

 The construction for $M_{\gl, c_1,c_2,c_3,W}$ is analogous, but we should also 
take into account the Heisenberg part. Recall that $V_{\gl,c_1,c_2}$ is itself
a tensor product of two VOAs:  $V_{\gl,c_1,c_2} = V_{\slN, c_1} \otimes 
V_{\h,c_2}$. The module $M_{\gl, c_1,c_2,c_3,W}$ that we are going to introduce
will be a tensor product of a module $M_{\slN,c_1,W}$ for $V_{\slN, c_1}$
and a module $M_{\h,c_2,c_3}$ for $V_{\h,c_2}$. The construction for
$M_{\slN,c_1,W}$ is completely parallel to the one we used above for $M_{\wdg,c,\dV}$.
Let us describe the Heisenberg VOA module $M_{\h,c_2,c_3}$.

 We take a one-dimensional space $\C 1$ and make it into a module for
$I\otimes \C[t_0] \oplus \C C_2$ by $I 1 = c_3 1, C_2 1 = c_2 1,
I\otimes t_0\C[t_0] 1 = 0$. The induced module
\begin{equation}
{\text {Ind}}_{I\otimes \C[t_0] \oplus \C C_2}^\h
\end{equation}
is an irreducible module over the Heisenberg algebra and will be denoted
by  $M_{\h,c_2,c_3}$, unless $c_2 = 0$, in which case we set 
$M_{\h,0,c_3}$ to be the one-dimensional factor $\C 1$ of the induced module.
The module $M_{\h,c_2,c_3}$ is an irreducible module for the Heisenberg VOA
$V_{\h,c_2}$ and
\begin{equation}
M_{\gl, c_1,c_2,c_3,W} = M_{\slN,c_1,W} \otimes M_{\h,c_2,c_3}
\end{equation}
is an irreducible VOA module for $V_{\gl,c_1,c_2}$. Finally,
\begin{equation}
M = M_{c,c_1,c_2,c_3,\dV,W} =  M_{\wdg,c,\dV} \otimes V_{\Lor}^+ 
\otimes M_{\gl, c_1,c_2,c_3,W}
\end{equation}
is an irreducible module for the toroidal VOA $V_{c,c_1,c_2}$.

Now we are ready to state the main result of this section:
\begin{theorem} \label{mod}
Let $c,c_1,c_2,c_3$ be complex numbers such that $c\neq 0$, $c\neq -h^\vee,
c_1 \neq -N$.
Let $\dV$ be an irreducible $\dg$-module and $W$ be an irreducible
$sl_N$-module.

(1) the module $M_{c,c_1,c_2,c_3,\dV,W}$ is a module for the toroidal
Lie algebra $\g_\tau$, $\tau = \mu \tau_1 + \nu \tau_2$, where
$\mu = \frac{1 - c_1}{c}, \nu = \frac{\frac{c_1}{N} - \frac{c_2}{N^2}}{c}$.

(2) $M_{c,c_1,c_2,c_3,\dV,W}$ is irreducible as a  $\g_\tau$-module.

(3) $M_{c,c_1,c_2,c_3,\dV,W} \cong U(c,\dV,W,c_3)$ as $\g_\tau$-modules.
\end{theorem}

{\sl Proof.}
The module $M = M_{c,c_1,c_2,c_3,\dV,W}$ is a module for the toroidal VOA
$V = V_{c,c_1,c_2}$. Under the correspondence \eqref{vk0}-\eqref{rank}
the commutator relations \eqref{g1g2}-\eqref{d0d0} of $\g$ hold in 
$V$. Thus, by  preservation of identities in VOA modules(see \cite {Li} Lemma 2.3.5),
the same relations hold in $M$, which proves that 
$M$ is a  $\g_\tau$-module.

 The module $M$ is an irreducible module for the toroidal
VOA $V$. In order to show that it is also irreducible as a 
$\g$-module, it is enough to show that $V$ is generated as a VOA
by the set 
\begin{align}
S = \bigg\{ g(-1) e^\ma, e^\ma, a_j(-1)e^\ma,& 
b_j(-1) e^\ma + \sum_{p=1}^N m_p E_{pj}(-1) e^\ma  \nonumber\\
& \bigg| g\in\dg, \m\in\Z^N,
 j=1,\ldots,N \bigg\}
\end{align}
which corresponds to \eqref{vk0}-\eqref{vg}. 
Once this is established, every $\g$-invariant subspace in a VOA module
is $V$-invariant, since, by the associativity property of VOAs every concatenation 
of n-th products may be expressed through the composition of operators from
$\g$. Thus  irreducibility with respect to $V$ will imply irreducibility with respect to $\g$.

 It is sufficient to show that the set $S$ generates each of the factors
$V_{\wdg,c} \otimes 1 \otimes 1$, $1\otimes V_{\Lor}^+ \otimes 1$ and
$1\otimes 1\otimes V_{\gl,c_1,c_2}$. To get the first factor, we note that
$g(-1) \in S$ for all $g\in\dg$. The VOA $V_{\wdg,c}$ is spanned by
the elements
\begin{equation}
g_1(-n_1) \ldots g_k (-n_k) 1, \quad g_i \in \dg, h_i \in \N.
\end{equation}
However, since
\begin{equation}
g_1(-n_1) \ldots g_k (-n_k) 1 = 
g_1(-1)_{(-n_1)} \left( \ldots \left( g_k(-1)_{(-n_k)} 1\right) \ldots \right),
\end{equation}
 we get that $S$ generates $V_{\wdg,c} \otimes 1 \otimes 1$.
Since $a_j(-1), b_j(-1), e^\ma \in S$, we also get 
$1\otimes V_{\Lor}^+ \otimes 1$ by a similar argument.
In particular we see that $b_j (-1) e^\ma$ is in the space generated by $S$
for any $\m\in\Z^N$. Choosing $\m$ with $m_p = \delta_{kp}, p=1,\ldots,N$,
we get that both $b_j(-1) e^{a_k}$ and $b_j(-1) e^{a_k} + E_{kj}(-1) e^{a_k}$
are in the space generated by $S$, and so is $E_{kj}(-1) e^{a_k}$. Finally
$(E_{kj}(-1) e^{a_k})_{(-1)} e^{-a_k} = E_{kj}(-1)$, and hence we can generate
$1\otimes1\otimes V_{\gl,c_1,c_2}$ as well. This completes the proof of part 
(2) of the theorem.

To prove (3) we note that $M_{c,c_1,c_2,c_3,\dV,W}$ has a natural
$\Z$-grading corresponding to the degree in the $t_0$-variable and moreover its zero
component $M_0 = \dV \otimes \C[\Lor^+] \otimes W$ is isomorphic to $T$.
One can easily see that the action of $\g_0$ on $M_0$ coincides with
\eqref{act1}-\eqref{act2}. Thus $M$ is a homomorphic image of $U$.
However, since the $\Z^{N+1}$-gradings on both modules are compatible,
and $U$ is graded simple, we get an isomorphism $M \cong U$.
This also proves that for $c\neq -h^\vee$, $\mu \neq N + 1$,
$U$  is a simple module and not just graded simple.

\end{document}